\documentclass[a4paper,11pt]{article}               
\usepackage{amssymb,amsfonts,amsmath}
\usepackage{graphicx,color,changebar}

\def\Re{{\mathfrak{Re}}}
\def\Im{{\mathfrak{Im}}}
\newcommand{\M}{{\mathcal M}}
\newcommand{\X}{{\mathcal X}}

\newcommand {\rank}{\mathop{\mathsf{rank}}\nolimits}

\newtheorem{theorem}{Theorem}
\newtheorem{remark}{Remark}
\newtheorem{lemma}{Lemma}
\newtheorem{definition}{Definition}

\newtheorem{corollary}{Corollary}

\newcommand {\mat}      [1] {\left[\begin{array}{#1}}
\newcommand {\rix}          {\end{array}\right]}

\newcommand{\Hn}{ {\mathbb{H}_n} }   
\newcommand{\Hnpsd}{ {\mathbb{H}_n^{\raisebox{0.15em}{{\fontsize{3}{2}\selectfont $\geq$}}}} } 
\newcommand{\Hmpsd}{ {\mathbb{H}_m^{\raisebox{0.15em}{{\fontsize{3}{2}\selectfont $\geq$}}}} } 
\newcommand{\Hnpd}{ {\mathbb{H}_n^{\raisebox{0.2em}{{\fontsize{3}{2}\selectfont $>$}}}} }  
\newcommand{\Hmpd}{ {\mathbb{H}_m^{\raisebox{0.2em}{{\fontsize{3}{2}\selectfont $>$}}}} }  

\newcommand{\XWpd}{ {\mathbb{X}^{\raisebox{0.2em}{{\fontsize{3}{2}\selectfont $>$}}}} }  
\newcommand{\XWpdpd}{ {\mathbb{X}^{\raisebox{0.2em}{{\fontsize{3}{2}\selectfont $\gg$}}}} }  

\title{Robust  port-Hamiltonian representations of passive systems}                                        

\author{Christopher Beattie\footnotemark[1] \, and Volker Mehrmann\footnotemark[2] \, and Paul Van Dooren\footnotemark[3]
}
\begin{document}

\maketitle

\begin{abstract}                          
We discuss the problem of robust representations of
stable and passive transfer functions in particular coordinate systems,
and focus in particular on the so-called port-Hamiltonian representations.
Such representations are typically far from unique and the degrees of freedom
are related to the solution set of the so-called Kalman-Yakubovich-Popov
linear matrix inequality (LMI). In this paper we analyze robustness measures for
the different possible representations and relate it to quality functions
defined in terms of the eigenvalues of the matrix associated with the LMI. In particular, we look at the analytic center of this LMI. From this, we then derive inequalities for the passivity radius
of the given model representation.
\end{abstract}
\renewcommand{\thefootnote}{\fnsymbol{footnote}}

\footnotetext[1]{
Department of Mathematics, Virginia Tech, Blacksburg, VA 24061, USA.
\texttt{beattie@vt.edu}. Supported by {\it Einstein Foundation Berlin},
through an Einstein Visiting Fellowship.}

\footnotetext[2]{
Institut f\"ur Mathematik MA 4-5, TU Berlin, Str. des 17. Juni 136,
D-10623 Berlin, FRG.
\texttt{mehrmann@math.tu-berlin.de}. Supported by {\it Einstein Foundation Berlin} via the Einstein Center ECMath and by Deutsche Forschungsgemeinschaft via Project A02 within CRC 910 `Control of self-organized nonlinear systems'}
\footnotetext[3]{
Department of Mathematical Engineering, UCL, Louvain-La-Neuve, Belgium.
\texttt{vandooren.p@gmail.com}.
and by {\it Deutsche Forschungsgemeinschaft},
through CRC 910 `Control of self-organized nonlinear systems'.
}
\section{Introduction}
We consider realizations of linear dynamical systems that are variously characterized as positive real, passive, or port-Hamiltonian.
We restrict ourselves to \emph{linear time-invariant systems} represented as
\begin{equation} \label{statespace}
 \begin{array}{rcl} \dot x & = & Ax + B u \quad \mbox{with}\ x(0)=0,\\
y&=& Cx+Du,
\end{array}
\end{equation}
where $u:\mathbb R\to\mathbb{C}^m$,   $x:\mathbb R\to \mathbb{C}^n$,  and  $y:\mathbb R\to\mathbb{C}^m$  are vector-valued functions denoting, respectively, the \emph{input}, \emph{state},
and \emph{output} of the system. The coefficient matrices $A\in \mathbb{C}^{n \times n}$,   $B\in \mathbb{C}^{n \times m}$, $C\in \mathbb{C}^{m \times n}$, and  $D\in \mathbb{C}^{m \times m}$ are constant.
Real and complex $n$-vectors ($n\times m$ matrices) are denoted by $\mathbb R^n$, $\mathbb C^{n}$
($\mathbb R^{n \times m}$, $\mathbb{C}^{n \times m}$), respectively.
We refer to (\ref{statespace}) concisely in terms of a four-tuple of matrices describing the realization $\M:=\left\{A,B,C,D\right\}$.

Our principal focus is on the structure of \emph{passive} systems and relationships with positive-real transfer functions and port-Hamiltonian system representations.
A nice introduction to passive systems can be found in the seminal papers of Willems (\cite{Wil71},  \cite{Wil72a}, \cite{Wil72b}), where a general notion of system passivity is introduced and linked to related system-theoretic notions such as positive realness and stability. Willems refers to the earlier works of Kalman \cite{Kal63}, Popov \cite{Pop73}, and Yakubovich \cite{Yak62}, where versions of what are now called Kalman-Yakubovich-Popov (KYP) conditions were derived.
Renewed interest in these ideas came from the study of port-Hamiltonian (pH) systems, which may be viewed as particular parameterizations of passive systems that arise from certain energy-based modeling frameworks (see e.g. \cite{Sch04}, \cite{Sch13}, \cite{SchJ14}). The KYP conditions lead to characterizations of system passivity through the solution set of an associated \emph{linear matrix inequality (LMI)}.  The convexity of this solution set has led to extensive use of
convex optimization techniques in systems and control theory (see e.g. \cite{BoyB90}, \cite{NesN94}).

The solution set of the KYP-LMI leads to a natural parametrization of families of pH realizations for a given passive system.
With this observation, it is not surprising that some pH realizations of a given system reflect well the underlying robustness of passivity to system perturbations and that some pH realizations will do this better than others.  Our main result shows that the analytic center of certain barrier functions associated with the KYP-LMI leads to favorable pH realizations in this sense; we derive computable bounds for the passivity 
radii for these realizations.

The paper is organized as follows. In Section~\ref{sec:prelim}, we recall the KYP conditions and link the solution set of the KYP-LMI to different system realizations that reflect passivity. In Section~\ref{sec:analytic} we review some basic concepts in convex analysis
 and introduce the concept of an \emph{analytic center} associated with a barrier function for the KYP-LMI.
In Section~\ref{sec:passrad} we define the passivity radius for a model realization $\M$, measuring its robustness to system perturbations that may lead to a loss of passivity.  We show that the analytic center of the KYP-LMI yields a model representation with good robustness properties. In Section~\ref{sec:criteria} we consider other measures that could also serve as criteria for robustness of model passivity.
In Section~\ref{sec:num} we illustrate our analytic results with a few numerical examples. We conclude with Section~\ref{sec:conclusion},
offering also some points for further research.

\section{Positive-realness, passivity, and
\hfill \break \mbox{port-Hamiltonian systems} }\label{sec:prelim}

We restrict ourselves to linear time-invariant systems as in \eqref{statespace} which are \emph{minimal}, that is, the pair $(A,B)$ is \emph{controllable}  (for all $s\in \mathbb C$, $\rank \mbox{\small $[\,s I-A \quad B\,]$} =n$ ), and the pair $(A,C)$ is \emph{observable} ($(A^\mathsf{H},C^\mathsf{H})$ is controllable). Here, the Hermitian (or conjugate) transpose (transpose) of a vector or matrix $V$ is denoted by
$V^{\mathsf{H}}$ ($V^{mathsf{T}}$).  We also assume that $\rank B=\rank C=m$.
We require input/output port dimensions to be equal  ($m$) and for convenience we assume that the system is initially in a quiescent state, $x(0)=0$. By applying the Laplace transform to \eqref{statespace} and eliminating the state, we obtain the  \emph{transfer function}
\begin{equation} \label{ABCD}
\mathcal T(s):=D+C(sI_n-A)^{-1}B, 
\end{equation}
mapping the Laplace transform of $u$ to the Laplace transform of $y$.
$I_n$ is the identity matrix in $\mathbb{C}^{n \times n}$
(subsequently, the subscript may be omitted if the dimension is clear from the context).
On the imaginary axis $\imath \mathbb R$, $\mathcal T(\imath\omega)$
describes the \emph{frequency response} of the system.

We denote the set of Hermitian matrices in $\mathbb{C}^{n \times n}$ by $\Hn$.
Positive definiteness (semidefiniteness) of  $A\in \Hn$ is denoted by $A>0$ ($A\geq 0$).
The set of all positive definite (positive semidefinite) matrices in $\Hn$ is denoted $\Hnpd$ ($\Hnpsd$).  The real and imaginary parts of a complex matrix, $Z$, are written as $\Re (Z)$ and $\Im (Z)$, respectively.

We proceed to review briefly some representations of linear systems associated with the notion of passivity.
\subsection{Positive-real systems}\label{sec:posreal}
Consider a system $\M $ as in (\ref{statespace}) and its transfer function $\mathcal T$ as in \eqref{ABCD}.
%
%
%
\begin{definition}\label{def:posreal}
A transfer function $\mathcal T(s)$ is {\em positive real}
if the matrix-valued rational function
\begin{equation}\label{defphi}
\Phi(s):= \mathcal T^{\mathsf{H}}(-s) + \mathcal T(s)
\end{equation}
is positive semidefinite for $s$ on the imaginary axis:
$$
\Phi(\imath\omega)\in \Hmpsd \quad \mbox{ for all }\omega\in \mathbb{R}.
$$
%
$\mathcal T(s)$ is \emph{strictly positive real} if $\Phi(\imath \omega)\in \Hmpd $ for all $\ \omega\in \mathbb{R}$.
\end{definition}
For any $X \in \Hn$, define the matrix function
\begin{eqnarray*} \label{prls}
W(X) &:=& \left[
\begin{array}{cc}
-X\,A - A^{\mathsf{H}}X & C^{\mathsf{H}} - X\,B \\
C- B^{\mathsf{H}}X & D+D^{\mathsf{H}}
\end{array}
\right]\\
& =& W(0)- \left[
\begin{array}{cc}
X\,A + A^{\mathsf{H}}X &  X\,B \\
B^{\mathsf{H}}X & 0
\end{array}
\right].
\end{eqnarray*}
 From \eqref{ABCD} and (\ref{defphi}), simple manipulations produce
\begin{eqnarray*}
\Phi\,(s) &= &
\left[ \begin{array}{cc} B^{\mathsf{H}}(-s\,I_n - A^{\mathsf{H}})^{-1} & I_m  \end{array} \right]
\, W(0) \left[ \begin{array}{c} (s\,I_n -A)^{-1}B \\ I_m \end{array} \right] \nonumber \\
 & =&
\left[ \begin{array}{cc} B^{\mathsf{H}}(-s\,I_n - A^{\mathsf{H}})^{-1} & I_m  \end{array} \right]
\, W(X) \left[ \begin{array}{c} (s\,I_n -A)^{-1}B \\ I_m \end{array} \right]. \label{popovs}
\end{eqnarray*}
Define the \emph{matrix pencils}:
\[
\mathcal L_{0}(s) =
s \left[ \begin{array}{cc|c} 0 & I_n & 0\\[1mm]
	-I_n & 0 & 0\\[1mm]
	\hline \rule{0mm}{1.1em}
        0 & 0 & 0 \end{array} \right]
	-\left[ \begin{array}{cc|c} 0 & A & B \\[1mm]
	A^{\mathsf{H}} & 0 & C^{\mathsf{H}} \\[1mm]
        \hline \rule{0mm}{1.1em}
         B^{\mathsf{H}} & C & D+D^{\mathsf{H}}  \end{array} \right],
\]
and, for any $X\in \Hn$,
\[
\mathcal L_{X}(s) =
\left[ \begin{array}{cc|c} I_n & 0 & 0 \\
-X & I_n & 0 \\ \hline 0 & 0 & I_m  \end{array} \right] \mathcal L_{0}(s)
\left[ \begin{array}{cc|c} I_n & -X & 0 \\
0 & I_n & 0 \\ \hline 0 & 0 & I_m  \end{array} \right],
\]
Observe for any $X\in \Hn$, $\mathcal L_{X}(s)$ and $\mathcal L_{0}(s)$ are equivalent pencils, since they are related
via a congruence transformation. Note also that $\Phi(s)$ is
the \emph{Schur complement} associated with the $(3,3)$ block of $\mathcal L_{0}(s)$
(and hence, also of $\mathcal L_{X}(s)$ for any $X\in \Hn$).

If $\mathcal T(s)$ is  positive real,
then it is known \cite{Wil71} that there exists $X\in \Hn$ such that the KYP-LMI holds, namely
\begin{equation} \label{KYP-LMI}
W(X) \geq 0,
\end{equation}
and so a factorization must exist:
\begin{equation} \label{lw}
W(X) =
\left[ \begin{array}{c} L^{\mathsf{H}} \\ M^{\mathsf{H}} \end{array} \right]\,
\left[ \begin{array}{cc} L &  M \end{array} \right]
\end{equation}
for $L\in \mathbb{C}^{r \times n}$ and $M \in \mathbb{C}^{r \times m} $,
 where $r=\rank W(X)$.
  Introducing $\mathcal G(s)= L(s\,I_n - A)^{-1}B + M,$~one may then define the \emph{spectral factorization} of $\Phi(s)$ as
\begin{equation} \label{spectral}
\Phi(s)= \mathcal G^{\mathsf{H}}(- s) \mathcal G(s).
\end{equation}
Define the solution set and subsets to the KYP-LMI (\ref{KYP-LMI}):
\begin{subequations}\label{LMIsolnsets}
\begin{align}
&{\mathbb X}:=\left\{ X\in \Hn \left|\  W(X) \geq 0 \right.\right\}, \label{XsolnWpsd} \\[1mm]
&\XWpd :=\left\{ X\in \Hn \left|   W(X) \geq 0,\ X >0 \right.\right\} = \Hnpd \cap {\mathbb X}, \label{XpdsolnWpsd} \\[1mm]
&\XWpdpd :=\left\{ X\in \Hn \left|   W(X) > 0,\ X >0 \right.\right\}, \label{XpdsolWpd}
\end{align}
\end{subequations}
For each $X\in{\mathbb X}$, there is a factorization of $W(X)$ of the form (\ref{lw}) leading to a spectral factorization (\ref{spectral}) of $\Phi(s)$. We are mainly interested in $\XWpd$ and $\XWpdpd$, which are
respectively,  the set of positive-definite solutions to the KYP-LMI (\ref{KYP-LMI}), and the subset of those solutions for which the KYP-LMI (\ref{KYP-LMI}) holds \emph{strictly}. 

An important subset of ${\mathbb X}$ are those solutions to (\ref{KYP-LMI}) for which the
rank $r$ of $W(X)$ is minimal ({i.e.}, for which $r=\rank\Phi(s)$). Let $S:= D+D^{\mathsf{H}}=\lim_{s\rightarrow\infty}\Phi(s)$.
If $S$ is nonsingular, then
the minimum rank solutions in $\XWpd$
are those for which $\rank W(X) = \rank S  = m$, which in turn is the case
if and only if the Schur complement of $S$ in $W(X)$ is zero.  This Schur
complement is associated with the \emph{algebraic Riccati equation (ARE)}:
\begin{multline}
\mathsf{Ricc}(X) := -XA-A^{\mathsf{H}}X  \\ -(C^{\mathsf{H}}-XB)S^{-1}(C-B^{\mathsf{H}}X)=0.\label{riccati}
\end{multline}
Solutions to (\ref{riccati}) produce directly a spectral factorization of $\Phi(s)$,
Indeed,  each solution $X$ of~\eqref{riccati} corresponds to a
\emph{Lagrangian invariant subspace} spanned by the columns of $U:=\mat{cc} I_n & -X^{\mathsf{T}} \rix^{\mathsf{T}} $
that remains invariant under the action of the Hamiltonian matrix
\begin{equation}\label{HamMatrix}
H:=\mat{cc} A-B S^{-1} C & - B S^{-1} B^{\mathsf{H}} \\
C^{\mathsf{H}} S^{-1} C & -(A-B S^{-1} C)^{\mathsf{H}} \rix.
\end{equation}
$U$ satisfies $HU=U A_F$ for a \emph{closed loop matrix} $A_F=A-BF$ with $F := S^{-1}(C-B^{\mathsf{H}}X)$ (see e.g., \cite{FreMX02}).
%
%
Each solution $X$ of~\eqref{riccati}  could also be associated with an \emph{extended Lagrangian invariant subspace}
for the pencil $\mathcal{L}_{0}(s)$ (see \cite{BenLMV15}), spanned by the columns of
$ \widehat{U}:=\mat{ccc} -X^{\mathsf{T}}
& I_n & -F^{\mathsf{T}}  \rix^{\mathsf{T}}$.
 In particular, $\widehat{U}$ satisfies
\[
\left[ \begin{array}{ccc} 0 & A & B \\
	A^{\mathsf{H}} & 0 & C^{\mathsf{H}} \\  B^{\mathsf{H}} & C & S  \end{array} \right] \widehat{U}
  =\left[ \begin{array}{ccc} 0 & I_n & 0\\
	-I_n & 0 & 0\\ 0 & 0 & 0 \end{array} \right] \widehat{U} A_F,
\]
see also {e.g.} \cite{IonOW98,Wil71}.   The condition that $S$ is invertible is equivalent
to the condition that the pencil $\mathcal{L}_{0}(s)$ has \emph{differentiation index one},
{i.e.}, all eigenvalues at $\infty$ are semi-simple, \cite{KunM06}.
If $\mathcal{L}_{0}(s)$ has no purely imaginary eigenvalues, then there are ${2n}\choose{n}$ solutions $X\in \Hn$ of (\ref{riccati}), each associated with an appropriate choice of a Lagrangian invariant subspace for $\mathcal{L}_{0}(s)$. Every choice leads to different spectra for the closed loop matrix, $A_F$ (see \cite{FreMX02,Wil71} for a parametrization of all possible Lagrangian subspaces).
Among the possible solutions of (\ref{riccati}) there are two extremal solutions, $X_-$ and $X_+$.  $X_-$ leads to a closed loop matrix, $A_F$, with spectra in the (open) left half-plane; $X_+$ leads to $A_F$ with spectra in the (open) right half-plane.  All solutions $X$ of (\ref{riccati}) are bracketed by $X_-$ and $X_+$:
\begin{equation} \label{Xbounded}
0\leq X_- \leq X \leq X_+
\end{equation}
and so, in this special case the set ${\mathbb X}$ is bounded, but it may be empty or the solution may be unique if $X_-=X_+$, see Section~\ref{sec:analytic}.
%
%

\subsection{Passive systems}\label{sec:DisSys}
%
\begin{definition}\label{def:passive}
A system $\M :=\left\{A,B,C,D\right\}$ is \emph{passive} if there exists a state-dependent
\emph{storage function}, $\mathcal H(x)\geq 0$, such that for any $\mu,t_0\in \mathbb R$ with $\mu>t_0$,
 the \emph{dissipation inequality} holds:
\begin{equation} \label{supply} \mathcal H(x(\mu))-\mathcal H(x(t_0)) \le \int_{t_0}^{\mu} \Re (y(t)^{\mathsf{H}}u(t)) \, dt
\end{equation}
%
If for all $\mu>t_0$, the inequality in \eqref{supply}
is strict then the system is \emph{strictly passive}.
\end{definition}
In the terminology of (\cite{Wil72a}),
$\Re (y(t)^{\mathsf{H}}u(t))$ is the \emph{supply rate} of the system.
A general theory of dissipative systems (of which passive systems are a special case) was developed in the seminal papers \cite{Wil71,Wil72a,Wil72b}, where links to earlier work by Kalman, Popov, and Yakubovich and the KYP-LMI (\ref{KYP-LMI}) are given.  Note that the original definition of passivity given by Willems was for real systems; we reformulate it here for complex systems.
\begin{theorem}[\cite{Wil71}]\label{LMIWil}
Suppose the system $\M$ of (\ref{statespace}) is minimal. Then the KYP-LMI (\ref{KYP-LMI}): $W(X)\geq 0$,
%
%
has a solution $X\in \Hnpd$ if and only if $\M$ is a passive system. If this is the case, then
\begin{itemize}
\item $\mathcal H(x):=\frac{1}{2}x^{\mathsf{H}}Xx$ defines a storage function associated with the supply rate $\Re(y^\mathsf{H}u)$ satisfying the dissipation inequality \eqref{supply};\\
\item there exist maximal and minimal solutions $X_- \leq X_+$  in $\Hnpd$ of \eqref{KYP-LMI},
such that for all solutions, $X$, of \eqref{KYP-LMI}:
\[
 0 < X_- \leq X \leq X_+.
\]
\end{itemize}
\end{theorem}
%
Recall that a matrix $A\in \mathbb{C}^{n \times n}$  is \emph{asymptotically stable} if all its eigenvalues are in the open left half plane and \emph{(Lyapunov) stable} if all eigenvalues are in the closed left half plane with any eigenvalues occurring on the imaginary axis being semisimple.
Theorem~\ref{LMIWil} asserts that if $X>0$  is a solution of $W(X)\geq 0$, then the system $\M$ of \eqref{statespace} is stable and if it satisfies $W(X)> 0$, then it is asymptotically stable, since $\mathcal H(x)$ is a Lyapunov function for $\M$ which is strict if $W(X)> 0$, (see {e.g.} \cite{LanT85}). Note, however, that  for (asymptotic) stability of $A$ it is sufficient if the $(1,1)$ block of $W(X)$ is (positive definite) positive semidefinite.
%
\begin{corollary}\label{cor:pr}
Consider a minimal system $\M $ as in (\ref{statespace}). $\M$ is passive if and only if it is positive real
and stable. It is strictly passive if and only if it is strictly positive real and asymptotically stable. In the latter case, $X_+ -X_- >0$ (\cite{Wil71}).
\end{corollary}
Note that minimality is not necessary for passivity.  For example, the system
$\dot{x}=-x,y = u$  is both stable and passive but not minimal.  In this case, the KYP-LMI (\ref{KYP-LMI}) is satisfied with any (scalar) $X>0$, the Hamiltonian may be defined as
$\mathcal H(x)=\frac{X}{2} x(t)^2$, and the dissipation inequality evidently holds since for $t_1\geq t_0$,
{\small
\vspace{-5mm}\begin{align*}
\mathcal{H}(x(t_1))-\mathcal{H}(x(t_0)) & = \frac{X}{2} (x(t_0) e^{-(t_1-t_0)})^2 -\frac{X}{2} (x(t_0))^2 \\
 = \frac{X}{2} (x(t_0))^2& (e^{-2(t_1-t_0)}-1) \leq 0 \leq \int_{t_0}^{t_1} y(t) \, u(t) \, dt
\end{align*} }
%
\subsection{Port-Hamiltonian systems}\label{sec:PH}
%
\begin{definition}\label{def:ph}
A linear time-invariant \emph{port-Hamiltonian (pH) system} is one for which the following realization is possible:
\begin{equation} \label{pH}
 \begin{array}{rcl} \dot x  & = & (J-R)Q x + (G-K) u,\\
y&=& (G+K)^{\mathsf{H}}Q x+Du,
\end{array}
\end{equation}
where $Q=Q^{\mathsf{H}} >0$, $J=-J^{\mathsf{H}}$, and 
\[
\left[ \begin{array}{lc} R & K \\ K^{\mathsf{H}} & \mathsf{sym}(D) \end{array} \right] \geq 0 \quad \mbox{with}
\quad \mathsf{sym}(D)=\frac12(D+D^{\mathsf{H}})
\]
%
%
\end{definition}
Port-Hamiltonian systems were introduced  in \cite{Sch04} as a tool for energy-based modeling.  An energy storage function
$\mathcal{H}(x)=\frac12x^{\mathsf{H}}Qx$ plays a central role and under the conditions given, the dissipation inequality (\ref{supply}) holds and so pH systems are always passive.  Conversely, any passive system may be represented as a pH system ({\ref{pH}), see e.g.,  \cite{BeaMX15_ppt}.
We briefly describe the construction of such a representation: Suppose the model $\M$ of \eqref{statespace} is minimal and passive and let $X=Q\in \XWpd$ be a solution of the KYP-LMI \eqref{KYP-LMI}.
For this $Q$, define $J:=\frac12 (AQ^{-1}- Q^{-1}A^{\mathsf{H}})$, $R:=-\frac 12 (AQ^{-1}+ Q^{-1}A^{\mathsf{H}})$,
$K:= \frac12\left(Q^{-1}C^{\mathsf{H}}-B\right)$, and $G:= \frac12\left(Q^{-1}C^{\mathsf{H}}+B\right)$.
Direct substitution shows that (\ref{statespace}) may be written in the form of (\ref{pH}), with $J=-J^{\mathsf{H}}$, and
\[
\left[ \begin{array}{lc} R & K \\ K^{\mathsf{H}} & \mathsf{sym}(D) \end{array} \right] =
\frac12 \left[ \begin{array}{lc} Q^{-1} & 0 \\ 0  & I \end{array} \right] \ W(Q)\
\left[ \begin{array}{lc} Q^{-1} & 0 \\ 0  & I \end{array} \right]  \geq 0.
\] 

Another possible representation of a passive system as a standard pH system can be obtained by using a symmetric factorization of a solution $X$ of 
\eqref{KYP-LMI}:  $X=T^{\mathsf{H}}T$ with $T\in \mathbb{C}^{n \times n}$ (e.g., the Hermitian square root of $X$ or the Cholesky factorization of $X$ are two possibilities).
One defines a state-space transformation, $x_T=Tx$, leading to an equivalent realization in $T$-coordinates:
\[
\{A_T,B_T,C_T,D_T\} := \{TAT^{-1}, TB, CT^{-1}, D \}
\]
The associated KYP-LMI \eqref{KYP-LMI} with respect to the new coordinate system can be written as
\begin{eqnarray*} \label{prls}
W_T(\hat{X}) &:=& \left[
\begin{array}{cc}
-\hat{X}\,A_T - A_T^{\mathsf{H}} \hat{X} & C_T^{\mathsf{H}} - \hat{X}\,B_T \\
C_T- B_T^{\mathsf{H}} \hat{X} & D+D^{\mathsf{H}}
\end{array}
\right] \geq 0,
\end{eqnarray*}
but since $X=T^{\mathsf{H}}T$ is a solution to the KYP-LMI \eqref{KYP-LMI} in the original state-space coordinates, we have $\hat{X}=I$ and
\begin{eqnarray}  \nonumber
W_T(I)&=&\left[ \begin{array}{cccc} T^{-H} & 0\\ 0 & I_m
\end{array}
\right] W(X)
\left[ \begin{array}{cccc} T^{-1} & 0\\ 0 & I_m
\end{array}
\right] \geq 0.  \label{PH}
\end{eqnarray}
We can then use the Hermitian and skew-Hermitian part of $A_T$ to obtain a pH representation in $T$-coordinates:
$J_T:=\frac12 (A_T-A_T^{\mathsf{H}})$, $R_T:=-\frac 12 (A_T+ A_T^{\mathsf{H}})$,
$K_T:= \frac12\left(C_T^{\mathsf{H}}-B_T\right)$, $G_T:= \frac12\left(C_T^{\mathsf{H}}+B_T\right)$, and $Q_T=I$, so that
\begin{equation} \label{pHalt}
 \begin{array}{rcl} \dot{x}_T  & = & (J_T-R_T)x_T + (G_T-K_T) u,\\
y&=& (G_T+K_T)^{\mathsf{H}} x_T+Du,
\end{array}
\end{equation}
is a valid pH representation in $T$ state-space coordinates.
%

We have briefly presented three closely related concepts for a minimal linear time-invariant system of the form \eqref{statespace}, positive realness, passivity, and that the system has pH structure. All three properties can be characterized algebraically via the solutions of linear matrix inequalities, invariant subspaces of special even pencils, or solutions of Riccati inequalities. However, there typically is a lot of freedom in the representation of such systems. This freedom, which results from particular choices of solutions to the KYP-LMI as well as  subsequent  state space transformations, may be used to make the representation more robust 
to perturbations. In many ways the pH representation seems to be the most robust representation \cite{MehMS16,MehMS17} and it also has many other advantages: it encodes the geometric and algebraic properties directly in the properties of the coefficients \cite{SchJ14}; it allows easy ways for structure preserving model reduction \cite{GugPBS12,PolV10}; it easily extends to descriptor systems \cite{BeaMXZ17_ppt,Sch13}; and it greatly simplifies optimization methods for computing stability and passivity radii \cite{GilMS18,GilS16,GilS17,OveV05}.

The remainder of this paper will deal with the question of how to make use of this freedom in the state space transformation to determine a 'good' or even `optimal' representation as a pH system. To do this we study in the next section the set of solutions of the KYP-LMI \eqref{KYP-LMI} and in particular its \emph{analytic center}.

\section{The analytic center of the solution set $\mathbb{X}^>$} \label{sec:analytic}
Solutions of the KYP-LMI  \eqref{KYP-LMI} and of linear matrix inequalities 
are usually obtained via optimization algorithms, see {e.g.} \cite{BoyEFB94}.  A common approach involves defining a
a \emph{barrier function} $b:\mathbb{C}^{n \times n} \to \mathbb C$ that is defined and finite throughout the interior of the constraint set becoming infinite as the boundary is approached, and
then using this function in the optimization scheme. The minimum of the barrier function
itself is of independent interest and is called the \emph{analytic center} of the constraint set \cite{GenNV99}.

We have seen in Section~\ref{sec:DisSys} that for a system $\M $ that is minimal and strictly passive there  exists a (possibly large) class of state space transformations that  transform the system to pH form. This class is characterized by the set $\XWpd$ of positive definite solutions of (the strict version) of the linear matrix inequality \eqref{KYP-LMI}.
If the set $\XWpd$ is non-empty and bounded, then the \emph{barrier function}
\begin{equation*}
b(X) := - \log \det W(X),
\end{equation*}
is bounded from below, but becomes infinitely large when $W(X)$ becomes singular.
The analytic center of the domain $\XWpd$ is
the minimizer of this barrier function. To characterize the analytic center, we analyze the \emph{interior}  of the set $\XWpd$
given by
%
\begin{eqnarray*}
\mathrm{Int} \, \XWpd &:=& \left\{ X\in \XWpd \; |
\mbox{ there exists } \delta>0 \mbox{ such that }\right .\\ && \left .X+\Delta_X \in \XWpd
 \mbox{for all}\;  \Delta_X \in \Hn \mbox{ with } \|\Delta_X\|\le \delta  \right\}.
\end{eqnarray*}
%
Here $\|\Delta_X\|$ is the spectral norm of $\Delta_X$ given by the maximal singular value of $\Delta_X$.
We compare $\mathrm{Int} \, \XWpd$ with the open set
\[
\XWpdpd = \left\{ X\in \XWpd \; | \; W(X)> 0 \right\}.
\]
%
Since $b(X)$ is finite for all points in $\XWpdpd$, there is an open neighborhood where it stays bounded, and this implies that $\XWpdpd \subseteq \mathrm{Int} \, \XWpd$.
The converse inclusion is not necessarily true. For example, consider a $2\times 2$ transfer function having the form, $\mathcal T(s)=\mathrm{diag} (t(s), 0 )$,  where $t(s)$ is a scalar-valued, strictly passive transfer function.  The LMI is rank deficient for all $X\in \Hn$ (hence $\XWpdpd=\emptyset$) but there is a relative interior, since $t(s)$ is strictly passive. The characterization when both sets are equal is given by the following theorem.
\begin{theorem}\label{thm:interior}
Suppose the system $\M$ of (\ref{statespace}) is passive and $\rank(B)=m$.  Then $\XWpdpd \equiv \mathrm{Int}\,\XWpd$.
\end{theorem}
{\bf Proof:} If $\XWpd=\emptyset$  then $\XWpdpd=\emptyset$  as well.
Otherwise, pick an $X\in \mathrm{Int}\,\XWpd$ and suppose
that $W(X)$ is  positive semidefinite and singular.  Then there exists a nontrivial $0\neq [z_1^\mathsf{T},z_2^\mathsf{T}]^\mathsf{T}\in\mathsf{Ker}\, W(X)$ and
an $\varepsilon>0$ sufficiently small so that
if $\Delta X \in \Hn$ with $\|\Delta X\|_F\leq \varepsilon$ then  $X+\Delta X\in \XWpd$.
Observe that for all such $\Delta X $, we have $W(X+\Delta X)=W(X)+\Gamma(\Delta X)\geq 0$, where
 $ \Gamma(\Delta X)=-\left[\begin{array}{cc} \Delta X A+ A^{\mathsf{H}} \Delta X &  \Delta X B \\
  B^{\mathsf{H}}\Delta X  & 0 \end{array} \right] $,
%
 and so
 %
\begin{equation}
0\leq \left[\begin{array}{c} z_1\\ z_2 \end{array} \right]^{\mathsf{H}} 
W(X+\Delta X)\ \left[\begin{array}{c} z_1\\ z_2 \end{array} \right]
 =
\left[\begin{array}{c} z_1\\ z_2 \end{array} \right]^{\mathsf{H}}
\Gamma(\Delta X)\ \left[\begin{array}{c} z_1\\z_2 \end{array} \right] \label{nonstrict}
 \end{equation}
 %
If there was a choice for $\Delta X \in \Hn$ with $\|\Delta X\|_F\leq \varepsilon$ producing strict inequality in \eqref{nonstrict},
then we would arrive at a contradiction, since the choice $-\Delta X$ satisfies the same requirements yet violates the inequality.  Thus, we must have equality in \eqref{nonstrict}
for all $\Delta X \in \Hn$,
which in turn implies
\[
W(X+\Delta X)\left[\begin{array}{c} z_1\\ z_2 \end{array} \right]= \Gamma(\Delta X)\left[\begin{array}{c} z_1\\ z_2 \end{array} \right]=0.
\]
This 
means that
 $B^{\mathsf{H}}\Delta X\, u=0$ for all $\Delta X \in \Hn$.
 If $z_1=0$, then we find that $\Delta X B\, z_2=0$ for all $\Delta X \in \Hn$ which in turn means $B z_2=0$ and so, in light of the initial hypothesis on $B$, means that $z_2=0$ which is a contradiction, and thus,
we must conclude that $W(X)$ is nonsingular after all, hence positive definite.
 To eliminate the last remaining case, suppose that $z_1\neq 0$.  Choosing first $\Delta X=I$, we find that $z_1\perp \mathsf{Ran}(B)$.
 Pick $0\neq b\in  \mathsf{Ran}(B)$ and define  $\Delta X =I-2ww^{\mathsf{H}}$ with
 $w=\frac{1}{\sqrt{2}}(\frac{z_1}{\|z_1\|}-\frac{b}{\|b\|})$.  Then
 $B^{\mathsf{H}}\Delta X\, z_1 = \frac{\|z_1\|}{\|b\|}B^{\mathsf{H}}b =0$
 which implies that $z_1=0$, and so, $z=0$,
 $W(X)>0$, and again the assertion holds.
\hfill $\Box$

To characterize when $\XWpdpd \equiv \mathrm{Int}\,\XWpd\neq \emptyset$ is complicated and there is some confusion in the literature, because several factors may influence the solvability of the KYP-LMI.
It is clear that $S=D+D^{H}$ must be positive definite for a solution to be in $\XWpdpd$, but clearly this is not sufficient as is demonstrated by the simple system $\dot x=u$, $y=x+du$ (with $d>0$) which is minimal and has $X=1\in  \XWpd$ as only solution of the KYP-LMI, so  $\mathrm{Int} \, \XWpd=\emptyset$.
In this case the associated pencil $\mathcal L_0$ has one eigenvalue $\infty$ and the  purely imaginary eigenvalues $\pm i  d$. It is passive, but not strictly passive, and stable (but not asymptotically stable), which is in contradiction to many statements in the literature, see e.~g.~\cite{Gri04}, where unfortunately no distinction between passivity and strict passivity is made. The system is, furthermore, port-Hamiltonian with $J=0$, $R=0$, $B,C=1$, $Q=1$, $P=0$ and $D=1$, and satisfies the dissipation inequality. An analogous example is obtained with $A= \left [ \begin{array}{cc} 0 & 1 \\ -1 & 0 \end{array} \right ]$, $B^{\mathsf{H}}=C=\left [ \begin{array}{cc} 1 & 0  \end{array} \right ]$, and $D=1$.
Then, $X=I_2$ is the unique positive definite solution, $\mathrm{Int} \, \XWpd=\emptyset$, and there are double eigenvalues of $\mathcal L_0$ at $\pm i$. The system is not asymptotically stable and not strictly passive, but stable, passive and pH. If in this example we choose $D=0$, then still $X=I$ is the unique positive definite solution of the KYP-LMI, but now $\mathcal L_0$ has the only purely imaginary eigenvalue $0$ and two Kronecker blocks of size $2$ for the eigenvalue $\infty$. In this case the Riccati equation \eqref{riccati} cannot be be formed and there does not exist a two-dimensional extended Lagrangian invariant space associated with the stable eigenvalues.

\begin{remark}{\rm
Note that the solutions $X_+$ and $X_-$ of the Riccati equation \eqref{riccati} yield singular $W(X_+)$ and $W(X_-)$, and are thus  on the boundary of $\XWpd$,
even though they are positive definite, as was pointed out in Theorem \ref{LMIWil}.
}\end{remark}

In the sequel, we  assume that $\XWpdpd\neq \emptyset$, so that the analytic center of $\XWpd$ is well-defined, see also \cite{NesN94}, and we can compute it as a candidate for a `good' solution to the LMI (\ref{KYP-LMI}) yielding a robust representation. This requires the solution of an optimization problem.
Keep in mind that we have the following the set inclusions
\[
\mathbb{X}_{W} \subset   \mathrm{Int}\XWpd=\XWpdpd  \subset \XWpd \subset   \Hnpd    \subset \Hn.
\]
%
For $X,Y \in \Hn $ we define the \emph{Frobenius inner product}
%
\[
 \langle X,Y\rangle  := \mathsf{trace}\!\left(\Re(Y)^T\Re(X)+\Im(Y)^T\Im(X)\right),
\]
which has the properties $\langle X,Y \rangle= \langle Y,X \rangle$, $\|X\|_F= \langle X,X \rangle^{\frac{1}{2}}$, and $YZ \rangle= \langle YX,Z \rangle = \langle XZ,Y \rangle$.
%
%
%
%
%

The \emph{gradient} of the barrier function $b(X)$ with respect to $W$ is given by
\[
\partial b(X) / \partial W =  -W(X)^{-1}.
\]
Using the chain rule and the Frobenius inner product, it follows from \cite{NesN94}
that $X\in \mathbb{C}^{n \times n}$ is an extremal point of $b(X)$ if and only if
\[
\langle \partial b(X) / \partial W,  \Delta W(X)[\Delta_X] \rangle = 0 \quad \mbox{for all} \; \Delta_X \in \Hn,
\]
where $\Delta W(X)[\Delta_X]$ is the incremental step in the direction $\Delta_X$ given by
\[
\Delta W(X)[\Delta_X] = -\left[ \begin{array}{cc}
A^{\mathsf{H}}\Delta_X+\Delta_X A & \Delta_X B \\ B^{\mathsf{H}}\Delta_X & 0 \end{array} \right].
\]
For an extremal point it is then necessary that
\begin{equation}\label{orth}
\langle W(X)^{-1} , \left[ \begin{array}{cc}
A^{\mathsf{H}}\Delta_X + \Delta_X A & \Delta_X B \\ B^{\mathsf{H}}\Delta_X & 0 \end{array} \right] \rangle \ =0
\end{equation}
for all $\Delta_X\in \Hn$.
Defining $F := S^{-1}(C-B^HX)$,
$P :=-A^HX-XA-F^{\mathsf{H}}SF$, and $A_F := A-BF$, it has been shown in \cite{GenNV99} that \eqref{orth} holds if and only if $P$ is invertible and
\begin{equation} \label{skew}
A_F^{\mathsf{H}}P +PA_F =0.
\end{equation}
Note that $P$ is nothing but the evaluation of the Riccati operator \eqref{riccati} at $X$, and that $A_F$ is the corresponding closed loop matrix. For the solutions of the Riccati equation we have
$P=\mathsf{Ricc}(X)=0$ (so $P$ is not invertible) and the corresponding closed loop matrix
has all its eigenvalues equal to an adequate   subset of the eigenvalues of the Hamiltonian matrix $H$ in (\ref{HamMatrix}). For an interior point of $\XWpd$
we have $P=\mathsf{Ricc}(X)> 0$, and hence $P$ has an invertible
square root $P^{\frac{1}{2}}\in \Hnpd$. Multiplying
(\ref{skew}) on both sides with  $P^{-\frac{1}{2}}$
we obtain that
\[
P^{-\frac{1}{2}}A_F^{\mathsf{H}}P^{\frac{1}{2}} + P^{\frac{1}{2}}A_FP^{-\frac{1}{2}}=0.
\]
Thus, $\hat A_F:=P^{\frac{1}{2}}A_F P^{-\frac{1}{2}}$ is skew-Hermitian and therefore  $\hat A_F$ as well as $A_F$ have all its eigenvalues on the imaginary axis. Hence the closed loop matrix $A_F$ of the analytic center has a spectrum that is also 'central' in a certain sense.

It is important to note that
\[
 \det (W(X)) = \det (\mathsf{Ricc}(X)) \det S, 
\]
which implies that we are also finding a stationary point of
$\det(\mathsf{Ricc}(X))$, since $S$ 
is constant and invertible. Since  $P\in \Hnpd$, we can rewrite the equations defining the analytic center of
$\XWpd$ as the solutions $X\in \Hn$, $P\in \Hnpd$, $F\in \mathbb C^{m,n}$ of the system of matrix equations
\begin{eqnarray} \nonumber 
S F  &=& C-B^{\mathsf{H}}X,\\
 P   &=& -A^{\mathsf{H}}X-XA-F^{\mathsf{H}}SF,  \label{FXP}\\ \nonumber 
 0&=& P(A-BF)+(A^{\mathsf{H}}-F^{\mathsf{H}}B^{\mathsf{H}})P\\
 &=&PA_F+A_F^{\mathsf{H}}P.\nonumber
\end{eqnarray}
%
System (\ref{FXP}) can be used to determine a solution via an iterative method, that uses a starting value $X_0$  to compute $P_0,F_0$ and then consecutively solutions $X_i$, $i=1,2,\ldots$ 
followed by computing a new $P_i$ and $F_i$.

\begin{remark}\label{rem:evp}{\rm
For given matrices $P,F$ the solution $X$ of \eqref{riccati} can be obtained via 
the invariant subspace
\[
\left[ \begin{array}{ccc} 0 & I_n & 0\\
	-I_n & 0 & 0\\ 0 & 0 & 0 \end{array} \right]\mat{c} -X \\ \phantom{-}I_n \\ -F \rix Z=
\left[ \begin{array}{ccc} 0 & A & B \\
	A^{\mathsf{H}} & -P & C^{\mathsf{H}} \\  B^{\mathsf{H}} & C & S \end{array} \right]\mat{c} -X \\ \phantom{-}I_n \\ -F \rix,
\]
Computing this subspace for $P=\mathsf{Ricc}(X)=0$ allows to compute the extremal solutions $X_+$ and $X_-$ of \eqref{riccati}, which then can be used to compute a starting point for an optimization scheme, see \cite{BanMVN17_ppt} for details.
}
\end{remark}

\section{The passivity radius}\label{sec:passrad}
Our goal to achieve `good' or even `optimal' pH representations of a passive system can be realized in different ways. A natural measure for optimality is a large \emph{passivity radius}
$\rho_{\M}$, which is the smallest perturbation (in an appropriate norm) to the coefficients of a model $\M $ that makes the system non-passive. Computational methods to determine $\rho_{\M}$ were introduced in \cite{OveV05}, while the converse question, what is the nearest passive system to a non-passive system has recently been discussed in \cite{GilMS18,GilS17}.

Once we have determined a solution $X\in \XWpd$ to the LMI~(\ref{KYP-LMI}), we can determine the representation \eqref{pH} 
as in Section~\ref{sec:PH} and the system is automatically passive (but not necessarily strictly passive). For each such representation we can determine the passivity radius and then choose the most robust solution $X\in \XWpd$ under perturbations by maximizing the passivity radius or by minimizing the condition number of $X^{\frac 12}$, which is the transformation matrix to pH form, see Section~\ref{sec:cond}.

\subsection{The $X$-passivity radius}

Alternatively,  for $X\in \mathrm{Int}\XWpd$ we can determine the smallest perturbation ${\Delta_\M}$ of the system matrices $A,B,C,D$ of the model $\M $ that leads to a loss of positive definiteness of $W(X)$, because then we are on the boundary of the set of passive systems. This is a very suitable approach for perturbation analysis, since for fixed $X$ the matrix
\begin{equation*} \label{hx}
W(X) = \left[
\begin{array}{cc}
0 & C^{\mathsf{H}} \\
C & D+D^{\mathsf{H}}
\end{array}
\right]
- \left[\begin{array}{cc}
A^{\mathsf{H}}X + X\,A & X\,B \\
B^{\mathsf{H}}X & 0
\end{array}
\right]
\end{equation*}
is linear in the unknowns $A,B,C,D$ and when we perturb the coefficients, then we preserve strict passivity as long as
\begin{eqnarray*}
&& W_{\Delta_\M} (X) := \left[
\begin{array}{cc}
 0 & (C+\Delta_C)^{\mathsf{H}} \\
(C+\Delta_C) & (D+\Delta_D)+(D+\Delta_D)^{\mathsf{H}}
\end{array}
\right]\\
&& \qquad - \left[\begin{array}{cc}
(A+\Delta_A)^{\mathsf{H}}X + X\,(A+\Delta_A) & X\,(B+\Delta_B) \\
(B+\Delta_B)^{\mathsf{H}}X & 0
\end{array}
\right]>0.
\end{eqnarray*}
Hence, given $X\in \mathrm{Int}\XWpd$, we can look for the smallest perturbation $\Delta_\M$ to the model $\M $ that makes $\det (W_{\Delta_\M}(X)=0$. To measure the size of the perturbation $\Delta_\M$ of a state space model $\M $ ,we use the Frobenius norm
\[
 \|\Delta_\M \| := \left \|\left[\begin{array}{ccc}
\Delta_A & \Delta_B \\
\Delta_C & \Delta_D
\end{array}\right] \right \|_F.
\]
%
%
Defining for $X\in \mathrm{Int}\XWpd$ the  \emph{$X$-passivity radius} as
\[
	\rho_\M(X):= \inf_{\Delta_\M\in \mathbb C^{n+m,n+m}}\left\{ \| \Delta_\M \| \; | \; \det W_{\Delta_\M}(X) = 0\right\}.
\]
%
%
Note that in order to compute $\rho_\M(X)$ for the model $\M $, we must first find a point $X\in \mathrm{Int}\XWpd$, since $W(X)$  must be positive definite to start with and also $X$ should be positive definite to obtain a state-space transformation to pH form.

We have the following relation between the $X$-passivity radius and the usual passivity radius.
\begin{lemma}\label{passbound}
Consider a given model $\M$ . Then the passivity radius is given by
\begin{eqnarray*}
\nonumber
	\rho_{\M}&=& \sup_{X\in \mathrm{Int}\XWpd}\inf_{\Delta_\M\in \mathbb C^{n+m,n+m}}\{\| \Delta_\M \| | \det W_{\Delta_\M}(X)=0\}\\ &=& \sup_{X\in \mathrm{Int}\XWpd} \rho_{\M}(X).\label{passive}
	\end{eqnarray*}
\end{lemma}
{\bf Proof:}
	If for any given $X \in \mathrm{Int} \XWpd$ we have that  $\| \Delta_\M \|< \rho_{\M}(X)$, then all systems $\M +\Delta_\M$ with $\| \Delta_\M \| < \rho_\M(X)$ are strictly passive.
	Therefore $\rho_\M \ge \sup_{\mathrm{Int}\XWpd} \rho_\M(X)$. Equality follows, since there exists a perturbation $\Delta_\M$ of norm $\rho_\M$ for which there does not exist a point $X\in \mathrm{Int}\XWpd$ for which $W_{\Delta_\M}(X)> 0$, hence, this system is not strictly passive anymore.
\hfill $\Box$

We derive an explicit formula for the $X$-passivity radius based on a one parameter  optimization problem. For this, we rewrite the condition $W_{\Delta_\M} (X)>0$ as
\begin{eqnarray} \nonumber
&& \left[\begin{array}{cc} -X & 0 \\ 0 & I_m \end{array}\right]
\left[\begin{array}{cc} A + \Delta_A & B+\Delta_B \\ C+\Delta_C & D+\Delta_D \end{array}\right] \\
&& \label{wdelta} \qquad +
\left[\begin{array}{cc} A^{\mathsf{H}}+\Delta_A^{\mathsf{H}} & C^{\mathsf{H}}+\Delta_C^{\mathsf{H}} \\ B^{\mathsf{H}}+\Delta_B^{\mathsf{H}} & D^{\mathsf{H}}+\Delta_D^{\mathsf{H}} \end{array}\right]
\left[\begin{array}{cc} -X & 0 \\ 0 & I_m \end{array}\right] >0.
\end{eqnarray}
Setting
\begin{equation} \label{defWhatX}
\hat W:=W(X), \; \hat X := \left[\begin{array}{cc} X & 0 \\ 0 & I_m \end{array}\right], \; \Delta_T := \left[\begin{array}{cc} -\Delta_A & -\Delta_B \\ \Delta_C & \Delta_D \end{array}\right],
\end{equation}
inequality (\ref{wdelta}) can be written as the LMI
\begin{equation} \label{WDelta}
W_{\Delta_\M}  = \hat W+\hat X \Delta_T + \Delta_T^{\mathsf{H}} \hat X >0
\end{equation}
as long as the system is still passive. In order to violate this condition, we need to find the smallest $\Delta_T$ such that $\det W_{\Delta_\M} =0$. Factoring out $\hat W^{-\frac{1}{2}}$ on both sides of (\ref{WDelta})  yields the characterization
\begin{eqnarray}
&&\det\left(I_{n+m} +
\hat W^{-\frac{1}{2}}\hat X \Delta_T \hat W^{-\frac{1}{2}}+
\hat W^{-\frac{1}{2}}\Delta_T^{\mathsf{H}} \hat X \hat W^{-\frac{1}{2}}\right) \nonumber \\
&&\det\left(I_{n+m} +
\left[ \begin{array}{cc} \hat W^{-\frac{1}{2}}\hat X & \hat W^{-\frac{1}{2}}\end{array}\right]
\left[ \begin{array}{cc} 0 & \Delta_T \\ \Delta_T^{\mathsf{H}} & 0 \end{array}\right]
\left[ \begin{array}{cc} \hat X\hat W^{-\frac{1}{2}} \\ \hat W^{-\frac{1}{2}}\end{array}\right] \right)
\nonumber  \\ \nonumber
&&\det\left(I_{2(n+m)} +
\left[ \begin{array}{cc} 0 & \Delta_T \\ \Delta_T^{\mathsf{H}} & 0 \end{array}\right]
\left[ \begin{array}{cc} \hat X \hat W^{-\frac{1}{2}} \\ \hat W^{-\frac{1}{2}}\end{array}\right]
\left[ \begin{array}{cc} \hat W^{-\frac{1}{2}}\hat X & \hat W^{-\frac{1}{2}}\end{array}\right] \right)\\
&& =0.\label{rhoX}
\end{eqnarray}
The minimal perturbation $\Delta_T$ for which this is the case was described in \cite{OveV05}
using the following theorem, which we have slightly modified in order to take into account the positive semi-definiteness of the considered matrix.
\begin{theorem} \label{thm:OveV} Consider the matrices $\hat X, \hat W$ in (\ref{defWhatX}) and the pointwise positive semidefinite matrix function
\begin{equation}\label{defmgamma}
M(\gamma):= \left[ \begin{array}{cc} \gamma \hat X \hat W^{-\frac{1}{2}} \\ \hat W^{-\frac{1}{2}} / \gamma\end{array}\right]
\left[ \begin{array}{cc} \gamma \hat W^{-\frac{1}{2}} \hat X & \hat W^{-\frac{1}{2}} / \gamma \end{array}\right]
\end{equation}
in the real parameter $\gamma$. Then the largest eigenvalue $\lambda_{\max}(M(\gamma))$ is a \emph{unimodal function} of $\gamma$ ({i.e.} it is first monotonically decreasing and then monotonically increasing with growing $\gamma$). At the minimizing value $\underline \gamma$,  $M(\underline{\gamma})$ has an eigenvector $z$, {i.e.}
\[
 M(\underline{\gamma}) z = \underline\lambda_{\max} z, \quad z:=\left[ \begin{array}{cc} u \\ v \end{array}\right],
 \]
where
$  \|u\|_2^2=\|v\|_2^2=\frac{1}{2}$.
The minimum norm perturbation $\Delta_T$ is of rank $1$ and is given by $\Delta_T=2uv^{\mathsf{H}}/\underline{\lambda}_{\max}$. It has norm $1/\underline{\lambda}_{\max}$
both in the spectral norm and in the Frobenius norm.
\end{theorem}
{\bf Proof}
The proof for a slightly different formulation was presented in \cite{OveV05}. Here we therefore just present the basic idea in our formulation. Let $\Gamma:=\mathrm{diag}(\gamma I_{m+n}, \frac{1}{\gamma} I_{m+n})$, then $M(\gamma)=\Gamma H_1 \Gamma$, while
\[
\Gamma^{-1}\left[ \begin{array}{cc} 0 & \Delta_T \\ \Delta_T^{\mathsf{H}} & 0 \end{array}\right]\Gamma^{-1}=\left[ \begin{array}{cc} 0 & \Delta_T \\ \Delta_T^{\mathsf{H}} & 0 \end{array}\right].
\]
Setting
\[
K(\gamma):=\left( I_{2(n+m)} + \left[ \begin{array}{cc} 0 & \Delta_T \\ \Delta_T^{\mathsf{H}} & 0 \end{array}\right] M(\gamma) \right),
\]
then $\det(K(\gamma))$ is independent of $\gamma$. But the vector $z$ is in the kernel of $K(\gamma)$, which implies that also $K(1)$ is singular. The value of the norm of the constructed
$\Delta_T$ follows from the fact that the subvectors $u$ and $v$ must have equal norm at the minimum.
\hfill $\Box$

A bound for $\underline{\lambda}_{\max}$ in Theorem~\ref{thm:OveV} is obtained by the following result.
\begin{corollary}\label{cor:lev} Consider the matrices $\hat X, \hat W$ in (\ref{defWhatX}) and the pointwise matrix function $M(\gamma)$ as in (\ref{defmgamma}). The largest eigenvalue of $M(\gamma)$ is also the largest eigenvalue of
\[
\gamma^2 \hat W^{-\frac{1}{2}} \hat X^2 \hat W^{-\frac{1}{2}} + \hat W^{-1}/\gamma^2.
\]
A simple upper bound for $\underline{\lambda}_{\max}$ is given by $\underline{\lambda}_{\max}\le \frac{2}{\alpha\beta}$ where $\alpha^2:=\lambda_{\min}(\hat W)$ and $\beta^2=\lambda_{\min}(\hat X^{-1}\hat W\hat X^{-1})$. The corresponding lower bound for $\| \Delta_T \|_F$ then becomes
\[
  \rho_\M(X) = \min_{\gamma} \| \Delta_T \|_F \ge \alpha\beta/2.
\]
\end{corollary}
{\bf Proof} Clearly $\|\hat W^{-1}\|_2\le \frac{1}{\alpha^2}$ and $\|\hat W^{-\frac{1}{2}} \hat X^2 \hat W^{-\frac{1}{2}} \|_2\le \frac{1}{\beta^2}$. So if we choose $\gamma^2=\frac{\beta}{\alpha}$
	then
\begin{eqnarray*}
&& \min_{\gamma} \|\gamma^2 \hat W^{-\frac{1}{2}} \hat X^2 \hat W^{-\frac{1}{2}} + \hat W^{-1}/\gamma^2\|\\
 && \qquad\le  \| (\beta/\alpha)\hat W^{-\frac{1}{2}} \hat X^2 \hat W^{-\frac{1}{2}} + (\alpha/\beta)\hat W^{-1}\|\\	 &&\qquad \le  \frac{2}{\alpha\beta} . \qquad \Box
\end{eqnarray*}

	We can construct a perturbation $\Delta_T=\epsilon (\alpha\beta)vu^{\mathsf{H}}$ of norm $|\epsilon|(\alpha\beta)$ which makes the matrix  $W_{\Delta_\M}$ singular and therefore gives an upper bound for $\rho_M(X)$. To compute this perturbation, let $u$, $v$ and $w$ be vectors of norm $1$, satisfying
$\hat W^{-\frac{1}{2}}u=u/\alpha$, $\hat W^{-\frac{1}{2}} \hat X v=w/\beta$, $\Delta_T=\epsilon(\alpha\beta)vu^{\mathsf{H}}$, and $\epsilon u^{\mathsf{H}}w=-|\epsilon u^{\mathsf{H}}w|$,
{i.e.},  $u$, $v$ and $w$ are the singular vectors associated with the largest singular values
$1/\alpha$ of	$\hat W^{-\frac{1}{2}}$ and $1/\beta$ of $\hat W^{-\frac{1}{2}}\hat X$. Inserting these values in (\ref{rhoX}), it follows that
\begin{eqnarray*}
&& \det\left(I_{n+m} +
	\left[ \begin{array}{cc} \hat W^{-\frac{1}{2}}\hat X & \hat W^{-\frac{1}{2}}\end{array}\right]
	\left[ \begin{array}{cc} 0 & \Delta_T \\ \Delta_T^{\mathsf{H}} & 0 \end{array}\right]
	\left[ \begin{array}{cc} \hat X\hat W^{-\frac{1}{2}} \\ \hat W^{-\frac{1}{2}}\end{array}\right] \right) \\
&& \qquad =\det\left(I_{n+m} +
	\left[ \begin{array}{cc} w & u \end{array}\right]
	\left[ \begin{array}{cc} 0 & \epsilon \\ \overline \epsilon & 0 \end{array}\right]
	\left[ \begin{array}{cc} w^{\mathsf{H}} \\ u^{\mathsf{H}} \end{array}\right] \right) \\
&& \qquad =
	\det\left(I_{2} + \left[ \begin{array}{cc} 0 & \epsilon \\ \overline \epsilon & 0 \end{array}\right]
	\left[ \begin{array}{cc} w^{\mathsf{H}} \\ u^{\mathsf{H}} \end{array}\right] \left[ \begin{array}{cc} w &  u \end{array}\right]\right).
\end{eqnarray*}
If we now choose the argument of the complex number $\epsilon$ such that $\epsilon u^{\mathsf{H}}w$ is real and negative and the amplitude of $\epsilon$ such that $1=|\epsilon u^{\mathsf{H}}w|+ |\epsilon|$, 	then
\[
\det \left(I_{2} + \left[ \begin{array}{cc} \epsilon u^{\mathsf{H}}w & \epsilon\\ \overline \epsilon & \overline \epsilon w^{\mathsf{H}}u \end{array}\right]\right) = (1-|\epsilon u^{\mathsf{H}}w|)^2-|\epsilon|^2=0.
\]
Since $0\le |u^{\mathsf{H}}w| \le 1$, we have that $\frac{1}{2} \le |\epsilon| \le 1$ and thus we have the interval $\alpha\beta/2 \le \rho_\M(X) \le |\epsilon| \alpha\beta$, $\frac{1}{2} \le |\epsilon| \le 1$ in which the  $X$-passivity radius is located.
If $u$ and $w$ are linearly dependent, then this interval shrinks to a point and the estimate is exact.
We summarize these observations in the following theorem.
\begin{theorem}\label{thm:Xpassivity}
Let $\M=\{A,B,C,D\}$ be a given model and let $X\in \mathrm{Int}\XWpd$. Then the $X$-passivity radius $\rho_\M(X)$ at this $X$ is bounded by
\[
 \alpha\beta/2 \le \rho_\M(X) \le  \alpha\beta/(1+|u^{\mathsf{H}}w|),
 \]
where
$\alpha^2:=\lambda_{\min}(\hat W)$, $\beta^2=\lambda_{\min}(\hat X^{-1}\hat W\hat X^{-1})$, $\hat W^{-\frac{1}{2}}u=u/\alpha$, $\hat W^{-\frac{1}{2}} \hat X v=w/\beta$.
Moreover, if $u$ and $w$ are linear dependent, then $\rho_\M(X)=\alpha\beta/2$.
\end{theorem}
If we use these bounds for the passivity radius in a pH system, we have the following corollary.
%
\begin{corollary}\label{cor:xeqI} If for a given system $\mathcal M$ we have that $X=I_n\in \mathrm{Int}\XWpd$   then 
the corresponding representation of the system is port-Hamiltonian, {i.e.}, it has the representation $\M:= \{J-R,G-K,G^{\mathsf{H}}+K^{\mathsf{H}},S+N\}$  and the X-passivity radius 
is given by
 \[
 \rho_\M(I)=\lambda_{\min}W(I)=\lambda_{\min}\left[\begin{array}{cc} R & K \\ K^{\mathsf{H}} & S	 \end{array}\right].
 \]
Moreover, if $X=I_n$ is the analytic center of $\mathrm{Int}\XWpd$, then
$\rho_\M(I)$ equals the passivity radius $\rho_\M$ of $\M$.
\end{corollary}
{\bf Proof}
    This follows directly from Theorem~\ref{thm:Xpassivity}, since then  $\alpha=\beta$ and we can choose $u=w$.
\hfill $\Box$
\begin{remark}\label{rem:alphabeta}{\rm Considering a pH representation,
the conditions $\hat W\geq \alpha^2 I_{n+m}$ and $\hat X^{-1} \hat W \hat X^{-1}\geq \beta^2 I_{n+m}$ yield the necessary (but not sufficient) condition for passivity that
\[
\left[ \begin{array}{cc} \hat W & \alpha \beta I_{n+m} \\
    	\alpha \beta I_{n+m} &  \hat X^{-1} \hat W \hat X^{-1}	
    	\end{array} \right] \geq 0,
\]  	
Using the square root $\hat T:= \hat X^{\frac12}$ of  $\hat X>0$ and a congruence transformation, one finds that this is equivalent to
\[
 \left[ \begin{array}{cc}  \hat T^{-1} \hat W  \hat T^{-1} & \alpha \beta I_{n+m} \\
    	\alpha \beta I_{n+m} &  \hat T^{-1} \hat W \hat T^{-1}	
    	\end{array} \right] \geq 0,
\]
which implies that $\hat T^{-1} \hat W  \hat T^{-1} \geq \alpha\beta I_{n+m}$.

Defining $\xi:=\lambda_{\min}(\hat T^{-1} \hat W  \hat T^{-1})$, we then also obtain the inequality $\xi \ge \alpha\beta$. Because of \eqref{PH}, this is also equal to
\[
 \xi =\lambda_{\min} \left[\begin{array}{cc}
R & K \\ K^{\mathsf{H}} & S
\end{array}\right]
\]
which suggests that pH representations
are likely to guarantee a good passivity margin. In order to compute the optimal product $\xi=\alpha\beta$, we could maximize $\xi$ under the constraint
\[
\hat W-\xi\hat X = \left[ \begin{array}{cc} - XA- A^{\mathsf{H}} X - \xi X & C^{\mathsf{H}}- XB \\ C-B^{\mathsf{H}} X & S-\xi I_m \end{array}\right]> 0.
\]
}
\end{remark}

\subsection{A port-Hamiltonian barrier function}

From our previous discussions it appears that if we want to make sure that a state-space representation has a large passivity radius, we should not maximize the determinant of $W(X)$, but maximize instead
\begin{equation} \label{tilde}
\det \left ( \left[ \begin{array}{cc} X^{-\frac12} & 0 \\ 0 & I_m \end{array}\right]  W(X)\left[ \begin{array}{cc} X^{-\frac12} & 0 \\ 0 & I_m \end{array}\right] \right )
\end{equation}
under the constraint $X>0$ so that its square root $T=X^\frac12$ exists.
Equivalently, if $X>0$, we can maximize  the determinant of
\begin{eqnarray*} \tilde W(X)&:=& W(X)\left[ \begin{array}{cc} X^{-1} & 0 \\ 0 & I_m \end{array}\right]\\
&=&  \left[ \begin{array}{cc} -XAX^{-1}-A^{\mathsf{H}} & C^{\mathsf{H}}-XB \\ CX^{-1}-B^{\mathsf{H}} & S \end{array}\right]
\end{eqnarray*}
which has the same eigenvalues and the same determinant, but is expressed in terms of the variable $X$.

With this modified barrier function $\tilde b(X):=- \log \det \tilde W(X)$ we obtain the following formulas for the gradient of the barrier with respect to $\tilde W$ and the incremental step of $\tilde W(X)$ in the direction $\Delta_X$.
\begin{eqnarray*}
\partial \tilde b(X/\partial \tilde W) &=& -\tilde W(X)^{-{\mathsf{H}}} = -W(X)^{-1}\left[ \begin{array}{cc} X & 0 \\ 0 & I_m \end{array}\right], \\
\Delta \tilde W(X)[\Delta_X] &=&  \left[ \begin{array}{cc}
XAX^{-1}\Delta_X  - \Delta_X A & -\Delta_X B \\ -CX^{-1}\Delta_X & 0 \end{array} \right]\left[ \begin{array}{cc} X^{-1} & 0 \\ 0 & I_m \end{array}\right].
\end{eqnarray*}
Using again the chain rule, the necessary condition for an extremal point is then that for all $\Delta_X \in \Hn$, $<\partial \tilde b(X) / \partial \tilde W,  \Delta \tilde W(X)[\Delta_X]> = 0$, or equivalently
\begin{equation} \label{orthtilde}
< W(X)^{-1} , \left[ \begin{array}{cc} XAX^{-1}\Delta_X -\Delta_X A  & -\Delta_X B \\ -CX^{-1}\Delta_X & 0
\end{array} \right] > \ =0 \; .
\end{equation}
Defining $P$ and $F$ as before, and using that
\[
 W(X)^{-1} = \left[ \begin{array}{cc} I_n & 0 \\ -F & I_m \end{array}\right]
 \left[ \begin{array}{cc} P^{-1} & 0 \\ 0 & S^{-1} \end{array}\right]
  \left[ \begin{array}{cc} I_n & -F^{\mathsf{H}} \\ 0 & I_m \end{array}\right],
\]
it then follows that \eqref{orthtilde} holds if and only if $P$ is invertible and for all $\Delta_X \in \Hn$
we have
\[
<  P^{-1} , (XAX^{-1}+F^{\mathsf{H}}CX^{-1})\Delta_X  - \Delta_X (A-BF) > \ =0,
\]
or equivalently
\[
< \Delta_X , P^{-1} (XAX^{-1}+F^{\mathsf{H}}CX^{-1}) - (A-BF)P^{-1}   > \ =0,
\]
which can be expressed as
\begin{eqnarray*}
&& P[(X^{-1}A^{\mathsf{H}}X+X^{-1}C^{\mathsf{H}}F)-(A-BF)]\\
&& \qquad + [(XAX^{-1}+F^HCX^{-1})-(A-BF)^{\mathsf{H}}]P=0.
\end{eqnarray*}
Note that if one performs the coordinate transformation
\[
\{A_T,B_T,C_T,D_T\} := \{TAT^{-1}, TB, CT^{-1}, D \}
\]
where $T^2=X$, then  $P_T=T^{-1}PT^{-1}$ and $F_T=FT^{-1}$, which yields
the equivalent condition
\begin{eqnarray*}
&& P_T[(A_T^{\mathsf{H}}-A_T)+(C^{\mathsf{H}}_T+B_T)F_T]\\
 && \qquad +[(A_T^{\mathsf{H}}-A_T)+(B_T+C_T^{\mathsf{H}})F_T]^{\mathsf{H}}P_T=0.
\end{eqnarray*}
Moreover, we have that
\[
 F_T = S^{-1}(C_T-B_T^{\mathsf{H}}), \quad P_T= -A_T-A_T^{\mathsf{H}}-F_T^{\mathsf{H}}SF_T.
\]
Thus, if we use a pH representation $\M=\{A_T,B_T,C_T,D_T\}=\{J-R,G-K,(G+K)^{\mathsf{H}},D\}$, then at the analytic center of the modified barrier function, we have
\begin{eqnarray*}
SF_T&=&K^{\mathsf{H}}, \ P_T= R -F_T^{\mathsf{H}}SF_T, \\
0&=&  P_T(J-GF_T)+(J-GF_T)^{\mathsf{H}}P_T,
\end{eqnarray*}
and of course $X_T=I$, which implies that the passivity radius is given by
\[
 \lambda_{\min} \left[ \begin{array}{cc} R & K \\ K^{\mathsf{H}} & S \end{array}\right].
\]
On the other hand, since
we have optimized the determinant of $\tilde W(X)$ which has the same determinant as \eqref{tilde}, it follows that
$$  \det \tilde W(X)= \det \left[ \begin{array}{cc} R & K \\ K^{\mathsf{H}} & S \end{array}\right]
$$
and we can expect to have obtained a nearly optimal passivity margin as well.

\section{Other Radii} \label{sec:criteria}
	
Even though in this paper we are focusing on the passivity, we point out that pH representations also have other properties that are important to consider. In this section, we consider two such properties.

\subsection{The Stability Radius}

If a positive definite solution of the LMI (\ref{KYP-LMI}) exists, then it follows from the positive definiteness of the  $(1,1)$ block that  the system is asymptotically stable. Hence we can employ the same technique that we have used for the characterization of the passivity radius to bound the  \emph{stability radius}, {i.e.}, the smallest perturbation $\Delta_A$ that makes the system loose its asymptotic stability.  Introducing the positive definite matrices
\begin{eqnarray}\nonumber
V(X)&:=& -XA-A^{\mathsf{H}}X , \\ V_{\Delta_A}(X)&:=& V(X) - X\Delta_A-\Delta_A^{\mathsf{H}}X,
\label{VDelta}
\end{eqnarray}
we define the \emph{$X$-stability radius} as the smallest $\Delta_A$, for which $V_{\Delta_A}(X)$ looses its positive definiteness, {i.e.} for
%
$X\in \X^{>0}$ with $V(X)> 0$, the  $X$-stability radius is defined as
	\[
	\rho_A(X):= \inf_{\Delta_A\in \mathbb{C}^{n \times n}}\left\{ \| \Delta_A \| \; | \; \det V_{\Delta_A}(X) = 0\right\}.
	\]
	%
%
Note that \eqref{VDelta} is defined similar to  \eqref{WDelta}, except for a sign change and therefore as for the passivity radius we obtain the bound
\[
\alpha\beta/2 \le \rho_A(X) \le \alpha\beta/(1+|u^{\mathsf{H}}w|),
\]
where $\alpha^2 = \lambda_{\min}(V)$, $\beta^2 = \lambda_{\min}(X^{-1}VX^{-1})$, and where $u$, $v$ and $w$ are normalized vectors satisfying
\[
V^{-\frac{1}{2}}u=u/\alpha, \quad V^{-\frac{1}{2}} X v=w/\beta.
\]
\begin{corollary}\label{cor:xeqIstab} If for a given system $\mathcal M$ we have that $X=I_n\in \mathrm{Int}\XWpd$   then at this point the corresponding representation of the system is port-Hamiltonian, {i.e.}, it has the representation $\M:= \{J-R,G-K,G^{\mathsf{H}}+K^{\mathsf{H}},D\}$  and the X-stability radius of this model is given by
\[
	\rho_A(I)=\lambda_{\min}V(I)=\lambda_{\min} (R).
\]
\end{corollary}
{\bf Proof}
	This follows directly from Theorem~\ref{thm:Xpassivity}, since then  $\alpha=\beta$ and we can choose $u=w$.
\hfill $\Box$

\begin{remark}\label{rem:stabrad}{\rm
		It follows from the conditions $V\geq \alpha^2 I_{n+m}$ and $X^{-1} V X^{-1}\geq \beta^2 I_{n+m}$ that a necessary (but not sufficient) condition for stability is given by $T^{-1}VT^{-1}\geq \alpha\beta I_n$, where $T=X^{\frac12}$.
	}
\end{remark}
Another robustness measure for the transformation $T$ is to require that the \emph{field of values} $\{ x^{\mathsf{H}}A_Tx | x\in \mathbb C^n\}$ of the transformed matrix $A_T$ is as far left as possible into the left half plane. In other words, we want to minimize the real part of the right most \emph{Rayleigh quotient} of $A_T$ given by
\[
\min_{T s.t. T^2\in \XWpd} \{ \max_{x\neq 0, x\in \mathbb{C}^n} \Re(\frac{x^{\mathsf{H}}A_Tx}{x^Hx}) \}.
\]
Writing $A_T=J_T-R_T$ with $J_T=-J_T^{\mathsf{H}}$ and $R_T=R_T^{\mathsf{H}}$  we clearly only need to $x^{\mathsf{H}}R_Tx$, since $\Re (x^{\mathsf{H}}J_Tx)=0$. In other words, we want to determine
\[
\min_{T\in \XWpd} \{ \max_{x\neq 0, x\in \mathbb{C}^n} (\frac{x^{\mathsf{H}}R_Tx}{x^{\mathsf{H}}x}) \},
\]
which amounts to maximizing the smallest eigenvalue of the $(1,1)$ block of the LMI~(\ref{KYP-LMI}).
It therefore is an alternative to maximize the determinant of $W_T(X)$, since this will tend to maximize all of its eigenvalues, including those of the principal submatrices.

We are not advocating here to use either of these two approaches to compute the stability radius of our system, since we know that it is given explicitly by the formula

\[ \rho_A =  \min_{\omega\in \mathbb{R} } \sigma_{\min}(A-\imath \omega I_n).
\]
We just want to stress here that using a pH realization based on the analytic center will also yield a robust stability margin.

\subsection{Well conditioned state-space transformations}\label{sec:cond}

Since for any solution $X\in \XWpd$, $T=X^{\frac 12}$  yields a state space transformation to pH form, we can also try to optimize the condition number of $T$ or directly the condition number of $X=T^2\in \XWpd$ within the set described by $ X_- \leq X \leq X_+ $.

Let us first consider the special case that $X_+$ and $X_-$ commute. In this  case there exists a unitary transformation $U$ that  simultaneously diagonalizes both $X_-$ and $X_+$, {i.e.}, $
U^{\mathsf{H}}X_-U=\mathrm{diag}\{\lambda^{(-)}_1, \ldots, \lambda^{(-)}_n\}$ and  $U^{\mathsf{H}}X_+U=\mathrm{diag}\{\lambda^{(+)}_1, \ldots, \lambda^{(+)}_n\}$.
Since $ X_- \leq X \leq X_+ $, it follows that each eigenvalue $\lambda_i$ of $X$, $i=1,\ldots,n$ must lie in the closed interval
$\lambda_i \in [\lambda^{(-)}_i , \lambda^{(+)}_i]$, and that these intervals are nonempty. If there exists a point $\lambda$
in the intersection of all these intervals, then $X=\lambda I_n$ is an optimal choice and it has condition number $\kappa(X)=1$. If not, then there are at least two non-intersecting intervals, which implies that
$\lambda^{(+)}_{\min} < \lambda^{(-)}_{\max}$ and hence that the closed interval $[\lambda^{(+)}_{\min} , \lambda^{(-)}_{\max}]$ must then be non-empty. Moreover, it must also intersect each of the intervals $[\lambda^{(-)}_i , \lambda^{(+)}_i]$ in at least one point.
Thus, if we choose for any $i=1,\ldots,n$
\[
 \lambda_i \in [\lambda^{(-)}_i , \lambda^{(+)}_i] \cap [\lambda^{(+)}_{\min} , \lambda^{(-)}_{\max}],
\]
then the resulting matrix will have optimal condition number $\kappa(X)=
\frac{\lambda^{(-)}_{\max}}{\lambda^{(+)}_{\min}}$, and hence $\kappa(T)=\sqrt{\frac{\lambda^{(-)}_{\max}}{\lambda^{(+)}_{\min}}}$.
The proof that this is optimal follows from the Loewner ordering of positive semidefinite matrices. The largest eigenvalue of $X$ must be larger or equal to $\lambda_{\max}^{(-)}$ and the smallest eigenvalue of $X$ must be smaller or equal to  $\lambda_{\min}^{(+)}$.

If  $X_+$ and $X_-$ do not commute, then there still exists a (non-unitary) congruence transformation $L$ which simultaneously diagonalizes both $D^{(-)}:=L^{\mathsf{H}}X_-L$ and $D^{(+)}:=L^{\mathsf{H}}X_+L$ but the resulting diagonal elements $d^{(-)}_i$
and $d^{(+)}_i$ are not the eigenvalues anymore. Nevertheless, the same construction holds for any matrix $X$, but we cannot prove optimality anymore. On the other hand, we can guarantee the bound
\[
\kappa(T)\le \max \left (\kappa(L), \kappa(L) \sqrt{\frac{d^{(-)}_{\max}}{d^{(+)}_{\min}}}\right ).
\]

\section{Numerical examples}\label{sec:num}
In this section we present a few numerical examples for realizations that are construct on the basis of the analytic center.

We first look at a real scalar transfer function of first degree ($m$=$n$=1) because in this case both analytic centers that we proposed earlier, are easy to compute analytically. The transfer functions of interest are given by
\begin{eqnarray*} T(s)&=&d + \frac{cb}{s-a}, \\
 \Phi(s)&=&\left[ \begin{array}{cc} b/(-s-a) & 1  \end{array} \right]
\left[
\begin{array}{cc}
0 & c \\
c & 2d
\end{array}
\right].\left[ \begin{array}{c} b/(s-a) \\ 1 \end{array} \right].
\end{eqnarray*}
If we assume that it the system is strictly passive, then
\[
W(x) = \left[\begin{array}{cc} -2ax & c-bx \\ c-bx & 2d \end{array}\right]
\]
must be positive definite for some value of $x$. This implies that $d>0$ and that
$\det W(x)=-4adx-(c-bx)^2=-b^2x^2-2(2ad-cb)x-c^2 $ is positive for some value of $x$.
This implies that the discriminant $(2ad-cb)^2-b^2c^2 = 4a^2d^2-4abcd = 4ad(ad-bc)$ must be positive. Since the system is also stable, we finally have the following necessary and sufficient conditions for strict passivity of a real first degree scalar function:
\[
  a<0, \quad d>0, \quad \det  \left[\begin{array}{cc}
a & b \\ c & d \end{array}\right] < 0.
\]
It is interesting to point out that the transfer function $\Phi(\jmath\omega)$ is a non-negative and unimodal function of $\omega$ with extrema at $0$ and $\infty$.
We thus can check strict passivity by verifying the positivity of $\Phi(\jmath\omega)$ at these two values, and the stability of $a$~:
\[ \Phi(\infty) = d >0 , \ \Phi(0)= \frac{2(ad-cb)}{a} >0, \ a<0.
\]
Since the determinant is quadratic in $x$, it is easy to determine the analytic center  $x_*$ of the linear matrix inequality $W(x)>0$ and the corresponding feedback and Riccati operator:
\begin{eqnarray*} x_*&=&\frac{c}{b}-2d\frac{a}{b^2}, \quad
f=\frac{a}{b}, \quad  p=2d\frac{a^2}{b^2}-2c\frac{a}{b},\\
W(x) &=& \left[\begin{array}{cc}
4d\frac{a^2}{b^2}-2c\frac{a}{b} & 2d\frac{a}{b} \\ 2d\frac{a}{b} & 2d \end{array}\right]
= \left[\begin{array}{cc} 1 & \frac{a}{b} \\ 0 & 1 \end{array}\right].
\left[\begin{array}{cc}p & 0 \\ 0  & 2d \end{array}\right]
\left[\begin{array}{cc} 1 & 0 \\ \frac{a}{b} & 1\end{array}\right],
\end{eqnarray*}
which implies $\det H(x)=2d\cdot p$ and the strict passivity condition
\[ a < 0, \quad r > 0 \quad \mathrm{and}  \quad p= \frac{2a}{b^2}(ad-bc) >0.
\]
The strict passivity is lost when either one of the following happens
\[ d+\delta_d=0, \quad a+\delta_a = 0, \quad \det \left[ \begin{array}{cc} a+\delta_a & b+\delta_b \\ c+\delta_c & d+\delta_d \end{array} \right]=0.
\]
Therefore, it follows that
\[
\rho = \min(d, a, \sigma_2 \left[ \begin{array}{ccc} a & b \\ c & d  \end{array} \right]) =  \sigma_2 \left[ \begin{array}{ccc} a & b \\ c & d  \end{array} \right].
\]
But at the analytic center $x_*=(2da-cb)/b^2$, we have
\[
  \det
2d\left[ \begin{array}{cc} 2\frac{a^2}{b^2}-\frac{ac}{bd} & \frac{a}{b} \\ \frac{a}{b} & 1 \end{array} \right]=2\frac{da}{b^2}(\frac{a}{b}-\frac{c}{d})
\]
which shows that the positivity at the analytic center yields the correct condition for strict passivity of the model.

If we use the port-Hamiltonian barrier function, we have to consider only $x>0$
since  $\det \tilde W(x)=\det W(x) /x$. One easily checks that the derivative of
$\det \tilde W(x)$ has a positive zero at $x_*=|\frac{c}{b}|$, which eventually yields a balanced realization $\M_T=\{a,\sqrt{|bc|}, bc/\sqrt{|bc|},d\}$, and an improved passivity radius
\[
 \sigma_2 \left[ \begin{array}{ccc} a & bc/\sqrt{|bc|} \\ \sqrt{|bc|} & d  \end{array} \right].
\]

As second test case we look at a random numerical  model  $\{A,B,C,D\}$ in pH form of state dimension $n=6$ and input/output dimension $m=3$  via
\[
  \left[ \begin{array}{ccc} R & K \\ K^{\mathsf{H}} & S  \end{array} \right] := M M^{\mathsf{H}} ,
\]
where $M$ is a random $(n+m)\times(n+m)$ matrix generated in MATLAB. From this we then identified the model $A:=-R/2, B:=-C^{\mathsf{H}}:=-K/2$ and $D:=S/2$.
This construction guarantees us that $X_0=I_n$ satisfies the LMI positivity constraint for the model $\M:=\{A,B,C,D\}$. We then used the Newton iteration developed in \cite{BanMVN17_ppt} to compute the analytic center $X_c$ of the LMI
\[ W(X):= \left[
\begin{array}{cc}
-X\,A - A^{\mathsf{H}}X & C^{\mathsf{H}} - X\,B \\
C- B^{\mathsf{H}}X & S
\end{array}
\right] >0,
\]

using the barrier function $b(X):=-\ln \det W(X)$. We then determined the quantities $\alpha^2 := \lambda_{\min}()\hat W)$, $\beta^2 :=  \lambda_{\min}(\hat X^{-1}\hat W \hat X^{-1})$, and $\xi := \lambda_{\min}(\hat X^{-\frac12}\hat W \hat X^{-\frac12})$,
where $\tilde W:=W(X_c)$ and $\tilde X := \mathrm{diag} \{X_c,I_m\}$. The constructed matrix
\[
\hat X^{-\frac12}\hat W \hat X^{-\frac12}=\left[ \begin{array}{ccc} R_c & K_c \\ K_c^{\mathsf{H}} & S_c  \end{array} \right]
\]
also contains the parameters of the port-Hamiltonian realization at the analytic center $X_c$. The results are given in the table below
\[
\begin{array}{c|c|c|c|c|c}
	 \alpha^2 & \beta^2 & \xi & \alpha\beta & \lambda_{\min}(R_c)  & \rho_{stab} \\
	 \hline
	 0.002366 & 0.001065 &  0.002381 & 0.001587 & 0.1254 & 0.1035
\end{array}
\]
They indicate that $\lambda_{\min}(R_c)$ at the analytic center is a good approximation of the true stability radius, and that $\lambda_{\min}(\hat X^{-\frac12}\hat W \hat X^{-\frac12})$ at the analytic center is a good approximation of the $X_c$-passivity radius estimate $\alpha\beta$.

\section{Conclusion} \label{sec:conclusion}

In this paper we have introduced the concept of the analytic center for a barrier function derived from the KYP LMI $W(X)>0$ for passive systems.
We have shown that the analytic center yields very promising results for choosing the  coordinate transformation  to pH form for a given passive system.

Several important issues have not been addressed yet. Can we also apply these ideas to the limiting situations where $W(X)$ is singular and/or $X$ is singular, or when the given system is not minimal.  More importantly, one should also analyze if these ideas can also be adapted to models
represented in descriptor form. 

Another interesting issue is that of finding the nearest passive system to a given system that is not passive. This has an important application in identification,
where may loose the property of passivity, due to computational and round-off errors incurred during the identification.
%
\section*{Acknowledgment}
The authors greatly appreciate the help of Daniel Bankmann from TU Berlin for the computation of the numerical example.

\bibliographystyle{plain}

\end{document}